\documentclass[12pt]{article}
\usepackage{amssymb}
\usepackage{mathrsfs}
\usepackage[hypertex]{hyperref}
\usepackage{amsfonts}
\usepackage{graphicx}
\usepackage{mathptmx}
\usepackage{latexsym,amsmath,amssymb,amsfonts,amsthm}

\newtheorem{theorem}{Theorem}[section]

\newtheorem{corollary}[theorem]{Corollary}

\newtheorem{remark}[theorem]{Remark}

\newtheorem{defn}[theorem]{Definition}

\begin{document}
\setcounter{page}{1}
\title{The Caffarelli-Kohn-Nirenberg inequalities and manifolds with nonnegative weighted Ricci curvature}
\author{Jing Mao}
%{\emph{This paper is dedicated to my family.}}
\date{}
\protect\footnotetext{\!\!\!\!\!\!\!\!\!\!\!\!{ MSC 2010: 53C21;
31C12.}
 }
\maketitle ~~~\\[-15mm]

\begin{center}
{\footnotesize  Department of Mathematics, Harbin Institute of
Technology (Weihai), Weihai, 264209, China \\
Email: jiner120@163.com, jiner120@tom.com }
\end{center}

%\\[5mm]

\begin{abstract}
We prove that $n$-dimensional ($n\geqslant3$) complete and
non-compact metric measure spaces with non-negative weighted Ricci
curvature in which some Caffarelli-Kohn-Nirenberg type inequality
holds are \emph{close} to the model metric measure $n$-space (i.e.,
the Euclidean metric $n$-space).
 \end{abstract}

\markright{\sl\hfill  J. Mao \hfill}

\section{Introduction}
\renewcommand{\thesection}{\arabic{section}}
\renewcommand{\theequation}{\thesection.\arabic{equation}}
\setcounter{equation}{0} \setcounter{maintheorem}{0}

Denote by $C_{0}^{\infty}(\mathbb{R}^{n})$ the space of smooth
functions with compact support in the $n$-dimensional Euclidean
space. Let $n\geqslant3$ be an integer and let $a$, $b$, $p$ be
constants satisfying the following conditions
\begin{eqnarray} \label{1.1}
-\infty<a<\frac{n-2}{2}, \qquad a\leqslant{b}\leqslant{a+1}, \qquad
p=\frac{2n}{n-2+2(b-a)}.
\end{eqnarray}
For all $u\in{C}_{0}^{\infty}(\mathbb{R}^{n})$, Caffarelli, Kohn and
Nirenberg \cite{ckn} has proven that there exists a positive
constant $C$ depending only on constants $a$, $b$ and $n$ (these
constants satisfy (\ref{1.1}) above) such that the functional
equality
\begin{eqnarray} \label{1.2}
\left(\int_{\mathbb{R}^{n}}|x|^{-bp}|u|^{p}dv_{\mathbb{R}^{n}}\right)^{\frac{1}{p}}\leqslant{C}\left(\int_{\mathbb{R}^{n}}|x|^{-2a}|\nabla{u}|dv_{\mathbb{R}^{n}}\right)^{\frac{1}{2}}
\end{eqnarray}
holds, where $|x|$ is the Euclidean length of $x\in\mathbb{R}^{n}$,
and $dv_{\mathbb{R}^{n}}$ is the Euclidean volume element determined
by the standard Euclidean metric. We know that when $a=b=0$, the
Caffarelli-Kohn-Nirenberg type inequality (\ref{1.2}) degenerates
into the classical Sobolev inequality; when $a=0$ and $b=1$, the
Caffarelli-Kohn-Nirenberg type inequality (\ref{1.2}) becomes the
Hardy inequality. The Sobolev and the Hardy inequalities have many
important applications (see, e.g.,
\cite{au2,au3,ckn,dcx,cw,hlp,h1,le,l,m3,m4} and the references
therein), so it is meaningful to investigate the
Caffarelli-Kohn-Nirenberg type inequality (\ref{1.2}). The sharpest
constant $C$ such that the inequality (\ref{1.2}) holds is called
the \emph{best} constant. In the study of functional inequalities,
finding the best constants is also interesting and difficult
subject.

Let $K_{a,b}$ be the best constant for the Caffarelli-Kohn-Nirenberg
type inequality (\ref{1.2}), which implies
\begin{eqnarray} \label{1.3}
K_{a,b}^{-1}=\inf\limits_{u\in{C}^{\infty}_{0}(\mathbb{R}^{n})-\{0\}}\frac{\left(\int_{\mathbb{R}^{n}}|x|^{-2a}|\nabla{u}|dv_{\mathbb{R}^{n}}\right)^{\frac{1}{2}}}{\left(\int_{\mathbb{R}^{n}}|x|^{-bp}|u|^{p}dv_{\mathbb{R}^{n}}\right)^{\frac{1}{p}}}.
\end{eqnarray}
There exist several conclusions related to the best constant
$K_{a,b}$. More precisely, for the Sobolev inequality (corresponding
to the case of $a=b=0$), Aubin \cite{au1} and Talenti \cite{t} have
\emph{separately} shown that
\begin{eqnarray*}
K_{0,0}=\left(\frac{1}{n(n-2)}\right)^{\frac{1}{2}}\left(\frac{\Gamma(n)}{nD_{n}\Gamma^{2}\left(\frac{n}{2}\right)}\right)^{\frac{1}{n}},
\end{eqnarray*}
where $D_{n}$ is the volume of the unit ball in $\mathbb{R}^{n}$,
and that a family of minimizers is given by
\begin{eqnarray*}
u(x)=\left(\lambda+|x|^{2}\right)^{1-\frac{n}{2}}, \qquad \lambda>0.
\end{eqnarray*}
For the case of $a=0$ and $0<b<1$, Lieb \cite{l} has proven that the
best constant $K_{0,b}$ is
\begin{eqnarray*}
K_{0,b}=\left(\frac{1}{(n-2)(n-bp)}\right)^{\frac{1}{2}}\left(\frac{(2-bp)\Gamma\left(\frac{2(n-bp)}{2-bp}\right)}{nD_{n}\Gamma^{2}\left(\frac{n-bp}{2-bp}\right)}\right)^{\frac{2(n-bp)}{2-bp}},
\end{eqnarray*}
and that a family of minimizers is given by
\begin{eqnarray*}
u(x)=\left(\lambda+|x|^{2-bp}\right)^{-\frac{n-2}{2-bp}}, \qquad
\lambda>0.
\end{eqnarray*}
Chou and Chu \cite{cc} have improved the above two cases to the
situation that $a\geqslant0$, $a\leqslant{b}<a+1$, and have shown
that the best constant $K_{a,b}$ is
\begin{eqnarray*}
K_{a,b}=\left(\frac{1}{(n-2a-2)(n-bp)}\right)^{\frac{1}{2}}\left(\frac{(2-bp+2a)\Gamma\left(\frac{2(n-bp)}{2-bp+2a}\right)}{nD_{n}\Gamma^{2}\left(\frac{n-bp}{2-bp+2a}\right)}\right)^{\frac{2(n-bp)}{2-bp+2a}},
\end{eqnarray*}
and that, for $a>1$, all minimizers are non-zero constant multiples
of the function
\begin{eqnarray*}
u(x)=\left(\lambda+|x|^{2-bp+2a}\right)^{-\frac{n-2-2a}{2-bp+2a}},
\qquad \lambda>0.
\end{eqnarray*}
Catrina and Wang \cite{cw} have investigated the remaining case of
the best constant $K_{a,b}$ and the existence or non-existence of
the minimizers.

From now on, we fix some notations, that is, let $n\geqslant3$ be an
integer, and let $a$, $b$ and $p$ be constants satisfying
\begin{eqnarray} \label{1.4}
0\leqslant{a}<\frac{n-2}{2}, \qquad a\leqslant{b}<{a+1}, \qquad
p=\frac{2n}{n-2+2(b-a)}.
\end{eqnarray}
For a prescribed complete manifold, denote by $C_{0}^{\infty}(M)$
the space of smooth functions with compact support on $M$, and let
$dv_{g}$ be the volume element (i.e., Riemannian measure) related to
the Riemannian metric $g$. In this paper, for convenience, \emph{we
make an agreement that $\mathrm{vol}(\cdot)$ represents the volume
of the given geometric object}.

Given a complete open manifold $M$ with non-negative Ricci
curvature, do Carmo and Xia \cite{dcx} have revealed the following
potential relation between a Caffarelli-Kohn-Nirenberg type
inequality (on $M$) of the form (\ref{1.2}) and the geometric
property related to the volume of a geodesic ball on $M$.

\begin{theorem} (\cite{dcx})  \label{theorem1}
Let $C_{1}\geqslant{K_{a,b}}$ be a constant, with $K_{a,b}$
determined by (\ref{1.3}) and (\ref{1.4}), and $M$ be an
$n$-dimensional ($n\geqslant3$) complete open manifold with
non-negative Ricci curvature. Fix a point $x_{0}\in{M}$ and denote
by $\rho$ the distance function on $M$ from $x_{0}$. Assume that,
for any $u\in{C_{0}^{\infty}(M)}$, we have
\begin{eqnarray*}
\left(\int_{M}t^{-bp}|u|^{p}dv_{g}\right)^{\frac{1}{p}}\leqslant{C_{1}}\left(\int_{M}t^{-2a}|\nabla{u}|^{2}dv_{g}\right)^{\frac{1}{2}}.
\end{eqnarray*}
Then for any $x\in{M}$, we have
\begin{eqnarray*}
\mathrm{vol}[B(x_{0},r)]\geqslant\left(\frac{K_{a,b}}{C_{1}}\right)^{\frac{n}{1+a-b}}\cdot{V}_{0}(r),
\quad \forall{r}>0,
\end{eqnarray*}
where $V_{0}(r)$ is the volume of an $r$-ball in $\mathbb{R}^{n}$.
\end{theorem}

For the special case that $a=b=0$, the above theorem is covered by
\cite[Theorem 2]{x2}. We prefer to point out one thing here, that
is, \cite[Theorem 2]{x2} has been improved by Mao \cite{m3} recently
(see \cite[Theorem 1.3]{m3} for the precise statement or the end of
Section 1 of \cite{m4} for the detailed explanation).

The purpose of this paper is to generalize Theorem \ref{theorem1}
above. For that, we need to use the following notions of smooth
metric measure spaces and the weighted Ricci curvature.

A smooth metric measure space (also known as the weighted measure
space) is actually a Riemannian manifold equipped with some measure
which is conformal to the usual Riemannian measure. More precisely,
for a given complete $n$-dimensional Riemannian manifold $(M,g)$
with the metric $g$, the triple $(M,g,e^{-f}dv_{g})$ is called a
smooth metric measure space, where $f$ is a \emph{smooth
real-valued} function on $M$ and, as before, $dv_{g}$ is the
Riemannian volume element related to $g$ (sometimes, we also call
$dv_{g}$ the volume density). Correspondingly, for a geodesic ball
$B(x_{0},r)$ on $M$, with center $x_{0}\in{M}$ and radius $r$, one
can also define its \emph{weighted (or $f$-)volume}
$\mathrm{vol}_{f}[B(x_{0},r)]$ as follows
\begin{eqnarray*}
\mathrm{vol}_{f}[B(x_{0},r)]:=\int\limits_{B(x_{0},r)}e^{-f}dv_{g}.
\end{eqnarray*}
Now, for convenience, \emph{we also make an agreement that in this
paper $\mathrm{vol}_{f}(\cdot)$ represents the weighted (or
$f$-)volume of the given geometric object on a metric measure
space}.

For a given smooth metric measure space $(M,g,e^{-f}dv_{g})$, the
following $N$-Bakry-\'{E}mery tensor
\begin{eqnarray*}
\mathrm{Ric}^{N}_{f}:=\mathrm{Ric}+\mathrm{Hess}f-\frac{df\otimes{df}}{N},
\end{eqnarray*}
with $\mathrm{Ric}$ and $\mathrm{Hess}$ the Ricci and the Hessian
operators on $M$, can be considered. Especially, when $N=\infty$,
the $N$-Bakry-\'{E}mery tensor $\mathrm{Ric}^{N}_{f}$ degenerates
into the so-called \emph{$\infty$-Bakry-\'{E}mery Ricci tensor}
$\mathrm{Ric}_{f}$ which is given by
\begin{eqnarray*}
\mathrm{Ric}_{f}=\mathrm{Ric}+\mathrm{Hess}f.
\end{eqnarray*}
The $\infty$-Bakry-\'{E}mery Ricci tensor is also called \emph{the
weighted Ricci tensor}. Bakry and \'{E}mery \cite{be1,be2}
introduced firstly and extensively investigated the generalized
Ricci tensor above and its relationship with diffusion processes.

Similar to the $p$-norm of smooth functions with compact support on
the manifold $(M,g)$, for the smooth metric measure space
$(M,g,e^{-f}dv_{g})$ and any $u\in{C_{0}^{\infty}(M)}$, we can
define \emph{the weighted $p$-norm} $\|u\|_{p;MMS}$ of $u$ as
follows
\begin{eqnarray*}
\|u\|_{p;MMS}:=\left(\int\limits_{M}|u|^{p}\cdot{e}^{-f}dv_{g}\right)^{\frac{1}{p}}.
\end{eqnarray*}
Clearly, when $f\equiv0$, the weighted $p$-norm is just the
$p$-norm.

Maybe people would have an illusion that smooth metric measure
spaces are not necessary to study since they are simply obtained
from correspondingly Riemannian manifolds by adding a conformal
measure to the Riemannian measure. However, the truth is not like
this, and they do have many differences. For instance, when
$\mathrm{Ric}_{f}$ is bounded from below, the Myer's theorem,
Bishop-Gromov's volume comparison, Cheeger-Gromoll's splitting
theorem and Abresch-Gromoll excess estimate cannot hold as the
Riemannian case. Here, for the purpose of comprehension, we would
like to repeat an example given in \cite[Example 2.1]{ww}. That is,
for the metric measure space
$(\mathbb{R}^{n},g_{\mathbb{R}^{n}},e^{-f}dv_{g_{\mathbb{R}^{n}}})$,
where $g_{\mathbb{R}^{n}}$ is the usual Euclidean metric and
$dv_{g_{\mathbb{R}^{n}}}$, as before, is the Euclidean volume
density related to $g_{\mathbb{R}^{n}}$, if
$f(x)=\frac{\lambda}{2}|x|^{2}$ for $x\in\mathbb{R}^{n}$, then we
have $\mathrm{Hess}=\lambda{g_{\mathbb{R}^{n}}}$ and
$\mathrm{Ric}_{f}=\lambda{g_{\mathbb{R}^{n}}}$. Therefore, from this
example, we know that unlike in the case of Ricci curvature bounded
from below uniformly by some positive constant, a metric measure
space is not necessarily compact provided
$\mathrm{Ric}_{f}\geqslant\lambda$ and $\lambda>0$. So, it is
meaningful to study the geometry of smooth metric measure spaces.
For the basic and necessary knowledge about the metric measure
spaces, we refer readers to the excellent work \cite{ww} of Wei and
Wylie. The subject on the metric measure space and the related
weighted Ricci tensor occurs naturally in many different subjects
and has many important applications (see, e.g., \cite{jl,gp,ww}).

\begin{theorem} \label{theorem3}
Let $C_{2}\geqslant{K}_{a,b}$ be a constant, where $K_{a,b}$ is
determined by (\ref{1.3}) and (\ref{1.4}). Assume that
$(M,g,e^{-f}dv_{g})$ is an $n$-dimensional ($n\geqslant3$) complete
and noncompact smooth metric measure space with non-negative
weighted Ricci curvature. For a point $x_{0}\in{M}$ at which
$f(x_{0})$ is away from $-\infty$, assume that the radial derivative
$\partial_{t}f$ satisfies $\partial_{t}f\geqslant0$ along all
minimal geodesic segments from $x_{0}$, with $t:=d(x_{0},\cdot)$ the
distance to $x_{0}$ (on $M$). If furthermore for any
$u\in{C_{0}^{\infty}(M)}$, the Caffarelli-Kohn-Nirenberg type
inequality
\begin{eqnarray}  \label{1.5}
\left(\int_{M}t^{-bp}|u|^{p}\cdot{e}^{-f}dv_{g}\right)^{\frac{1}{p}}\leqslant{C_{2}}\left(\int_{M}t^{-2a}|\nabla{u}|^{2}\cdot{e}^{-f}dv_{g}\right)^{\frac{1}{2}}
\end{eqnarray}
holds, then we have
\begin{eqnarray*}
\mathrm{vol}_{f}[B(x_{0},r)]\geqslant\left(\frac{K_{a,b}}{C_{2}}\right)^{\frac{n}{1+a-b}}\cdot{e}^{-f(x_{0})}\cdot{V}_{0}(r),
\quad \forall{r}>0,
\end{eqnarray*}
where $V_{0}(r)$ is the volume of an $r$-ball in $\mathbb{R}^{n}$.
\end{theorem}

\begin{remark}
\rm{ If $f\equiv0$ on $M$, then the metric measure space
$(M,g,e^{-f}dv_{g})$ can be seen as the Riemannian manifold $(M,g)$
directly. Clearly, in this case, Theorem \ref{theorem3} is totally
the same with Theorem \ref{theorem1} above. So, we can equivalently
say that \cite[Theorem 1.1]{dcx} is only a special case
 of Theorem \ref{theorem3}. Moreover, as pointed out in \cite[Remark 1.4]{m3} or
\cite[Remark1.4]{m4}, since $f$ is a smooth real-valued function on
the complete non-compact manifold $M$, we know that if $f(x)$ does
not tend to $-\infty$ as $x$ tends
 to the infinity, then $x_{0}$ can be chosen arbitrarily; if $f(x)\rightarrow-\infty$ as $x$ tends
 to the infinity, then $x_{0}$ can be chosen to any point
 except those points near the infinity. Besides, we say that Theorem
 \ref{theorem3} here is \emph{sharper} than Theorem \ref{theorem1},
 since at $x_{0}$, one can always set up a global polar
 coordinate chart $\{t,\xi\}$ with
 $(t,\xi)\in[0,+\infty)\times\mathbb{S}^{n-1}$ for the complete
 non-compact manifold $M$, where $t:=d(x_{0},\cdot)$ as in Theorem
\ref{theorem3}, and then we have
$e^{-f(t,\xi)}\leqslant{e}^{-f(0,\xi)}=e^{-f(x_{0})}$ by applying
the assumption $\partial_{t}f\geqslant0$ along all minimal geodesic
segments from $x_{0}$, which leads to
\begin{eqnarray*}
\mathrm{vol}[B(x_{0},r)]\geqslant{e^{f(x_{0})}}\cdot\mathrm{vol}_{f}[B(x_{0},r)]=e^{f(x_{0})}\int\limits_{B(x_{0},r)}e^{-f}dv_{g}\geqslant\left(\frac{K_{a,b}}{C_{2}}\right)^{\frac{n}{1+a-b}}\cdot{V}_{0}(r),
\quad \forall{r}>0
\end{eqnarray*}
by using the conclusion of Theorem \ref{theorem3} directly. The
assumption on $f$ (i.e., \emph{finite} at $x_{0}$ and \emph{monotone
non-decreasing} in the radial direction) in Theorem \ref{theorem3}
seems a little strong here. However, it is not difficult to see that
there are many examples satisfying this condition. For instance, as
mentioned in \cite[Remark 1.4]{m3}, in the polar coordinate chart
$\{t,\xi\}$ constructed above, one can choose $f(t,\xi)=t$ for
$t\geqslant0$. Even more, one can find that functions
$f(t,\xi)=t^{\ell}$, $\ell>1$, for $t\geqslant0$ are also
acceptable. Hence, from this aspect, Theorem \ref{theorem3}
generalizes Theorem \ref{theorem1} a lot.}
\end{remark}

By applying Theorem \ref{theorem3} above, \cite[Theorem 3.3,
Corollary 3.4 and Theorem 4.2]{fmi} (see also Theorems
\ref{theorem4} and \ref{theorem5} in Section 2), and \cite[Theorem
1.2]{ww} (see also Theorem \ref{theorem6} in Section 2), we can
prove the following rigidity theorem.

\begin{corollary} \label{corollarys1}
Assume that $(M,g,e^{-f}dv_{g})$ is an $n$-dimensional
($n\geqslant3$) complete and noncompact smooth metric measure space
with non-negative weighted Ricci curvature. For a point
$x_{0}\in{M}$ at which $f(x_{0})$ is away from $-\infty$, assume
that the radial derivative $\partial_{t}f$ satisfies
$\partial_{t}f\geqslant0$ along all minimal geodesic segments from
$x_{0}$, with $t:=d(x_{0},\cdot)$ the distance to $x_{0}$ (on $M$).
If furthermore for any $u\in{C_{0}^{\infty}(M)}$, the
Caffarelli-Kohn-Nirenberg type inequality
\begin{eqnarray*}
\left(\int_{M}t^{-bp}|u|^{p}\cdot{e}^{-f}dv_{g}\right)^{\frac{1}{p}}\leqslant{K_{a,b}}\left(\int_{M}t^{-2a}|\nabla{u}|^{2}\cdot{e}^{-f}dv_{g}\right)^{\frac{1}{2}}
\end{eqnarray*}
holds, where $K_{a,b}$ is determined by (\ref{1.3}) and (\ref{1.4}),
then  $(M,g)$ is isometric to
$\left(\mathbb{R}^{n},g_{\mathbb{R}^{n}}\right)$. Moreover, in this
case, we have $f\equiv{f}(x_{0})$ is a constant function with
respect to the variable $t$, and
$e^{-f}dv_{g}=e^{-f(x_{0})}dv_{g_{\mathbb{R}^{n}}}$. Here, as
before, $g_{\mathbb{R}^{n}}$ and $dv_{g_{\mathbb{R}^{n}}}$ are the
usual Euclidean metric and the Euclidean volume density related to
$g_{\mathbb{R}^{n}}$, respectively.
\end{corollary}

It is interesting to know \emph{under what kind of conditions} a
complete open $n$-manifold ($n\geqslant2$) is isometric to
$\mathbb{R}^{n}$ or has finite topological type, which in essence
has relation with the splittingness of the prescribed manifold. This
is a classical topic in the global geometry and has been studied
extensively (see, e.g., \cite{dcx2,m,pp}).

\section{Useful facts}
\renewcommand{\thesection}{\arabic{section}}
\renewcommand{\theequation}{\thesection.\arabic{equation}}
\setcounter{equation}{0} \setcounter{maintheorem}{0}

We would like to review \cite[Theorem 1.2]{ww}, which is the
\emph{cornerstone} of the proof of Theorem \ref{theorem3} shown in
the next section, and \cite[Theorem 3.3, Corollary 3.4 and Theorem
4.2]{fmi}, which are necessary to prove Corollary \ref{corollarys1}.
However, before that, some necessary preliminaries should be
introduced first. In fact, one can find more detailed versions (cf.
\cite[Section 2]{fmi}, \cite[Section 2]{m1} and \cite[Section 2.1 of
Chapter 2]{m2}) of the following preliminaries, but we still give a
simple version here so that readers can understand \cite[Theorem
1.2]{ww} and \cite[Theorem 3.3, Corollary 3.4 and Theorem 4.2]{fmi})
completely and clearly.

\subsection{Preliminaries}

Denote by $\mathbb{S}^{n-1}$ the unit sphere in $\mathbb{R}^{n}$.
Given an $n$-dimensional ($n\geqslant2$) complete Riemannian
manifold $(M,g)$ with the metric $g$, for a point $x\in{M}$, let
$S_{x}^{n-1}$ be the unit sphere with center $x$ in the tangent
space $T_{x}M$, and let $Cut(x)$ be the cut-locus of $x$, which is a
closed set of zero $n$-Hausdorff measure. Clearly,
\begin{eqnarray*}
\mathbb{D}_{x}=\left\{t\xi|0\leqslant{t}<d_{\xi},\xi\in{S_{x}^{n-1}}\right\}
\end{eqnarray*}
is a star-shaped open set of $T_{x}M$, and through which the
exponential map
$\exp_{x}:\mathbb{D}_{x}\rightarrow{M}\backslash{Cut(x)}$ gives a
diffeomorphism from $\mathbb{D}_{x}$ to the open set
$M\backslash{Cut(x)}$, where $d_{\xi}$ is defined by
\begin{eqnarray*}
d_{\xi}=d_{\xi}(x) : = \sup\{t>0|~\gamma_{\xi}(s):=
\exp_x(s\xi)~{\rm{is~ the~ unique}}~{\rm{minimal
~geodesic~joining}}~x ~ {\rm{and}}~\gamma_{\xi}(t)\}.
\end{eqnarray*}

As in \cite{ic}, we can introduce two
 important maps used to construct the geodesic spherical coordinate chart at a prescribed point
 on a Riemannian manifold. For a fixed vector $\xi\in{T_{x}M}$,
$|\xi|=1$, let $\xi^{\bot}$ be the orthogonal complement of
$\{\mathbb{R}\xi\}$ in $T_{x}M$, and let
$\tau_{t}:T_{x}M\rightarrow{T_{\exp_{x}(t\xi)}M}$ be the parallel
translation along $\gamma_{\xi}(t)$. The path of linear
transformations
$\mathbb{A}(t,\xi):\xi^{\bot}\rightarrow{\xi^{\bot}}$ is defined by
 \begin{eqnarray*}
\mathbb{A}(t,\xi)\eta=(\tau_{t})^{-1}Y_{\eta}(t),
 \end{eqnarray*}
where $Y_{\eta}(t)=d(\exp_x)_{(t\xi)}(t\eta)$ is the Jacobi field
along $\gamma_{\xi}(t)$ satisfying $Y_{\eta}(0)=0$, and
$(\nabla_{t}Y_{\eta})(0)=\eta$. Moreover, for $\eta\in{\xi^{\bot}}$,
set
$
\mathcal{R}(t)\eta=(\tau_{t})^{-1}R(\gamma'_{\xi}(t),\tau_{t}\eta)
\gamma'_{\xi}(t),
$
where the curvature tensor $R(X,Y)Z$ is defined by
$R(X,Y)Z=-[\nabla_{X},$ $ \nabla_{Y}]Z+ \nabla_{[X,Y]}Z$. Then
$\mathcal{R}(t)$ is a self-adjoint operator on $\xi^{\bot}$, whose
trace is the radial Ricci tensor
$\mathrm{Ric}_{\gamma_{\xi}(t)}\left(\gamma'_{\xi}(t),\gamma'_{\xi}(t)\right)$.
 Clearly, the map $\mathbb{A}(t,\xi)$ satisfies the Jacobi
equation
 $\mathbb{A}''+\mathcal{R}\mathbb{A}=0$ with initial conditions
 $\mathbb{A}(0,\xi)=0$, $\mathbb{A}'(0,\xi)=I$.
By Gauss's lemma, the Riemannian metric of $M \backslash Cut(x)$ in
the geodesic spherical coordinate chart can be expressed by
 \begin{eqnarray} \label{2.1}
 ds^{2}(\exp_{x}(t\xi))=dt^{2}+|\mathbb{A}(t,\xi)d\xi|^{2}, \qquad
\forall t\xi\in\mathbb{D}_{x}.
 \end{eqnarray}
We consider the metric components  $g_{ij}(t,\xi)$, $i,j\geq 1$, in
a coordinate system $\{t, \xi_a\}$ formed by fixing  an orthonormal
basis $\{\eta_a, a\geq 2\}$ of
 $\xi^{\bot}=T_{\xi}S^{n-1}_x$, and then extending it to a local frame $\{\xi_a, a\geq2\}$ of
$S_x^{n-1}$. Define
 a function $J>0$ on $\mathbb{D}_{x}\backslash\{x\}$ by
\begin{equation} \label{2.2}
J^{n-1}=\sqrt{|g|}:=\sqrt{\det[g_{ij}]}.
\end{equation}
 Since
 $\tau_t: S_x^{n-1}\to S_{\gamma_{\xi}(t)}^{n-1}$ is an
isometry, we have
$$
\langle d(\exp_x)_{t\xi}(t\eta_{a}),
d(\exp_x)_{t\xi}(t\eta_{b})\rangle_{g} =\langle
\mathbb{A}(t,\xi)(\eta_{a}), \mathbb{A}(t,\xi)(\eta_{b})\rangle_{g},
$$
and then $ \sqrt{|g|}=\det\mathbb{A}(t,\xi).$ So, by applying
(\ref{2.1}) and (\ref{2.2}), the volume $\mathrm{vol}(B(x,r))$ of a
geodesic ball $B(x,r)$, with radius $r$ and center $x$, on $M$ is
given by
\begin{eqnarray}\label{2.3}
\mathrm{vol}(B(x,r))=\int\limits_{S_{x}^{n-1}}\int\limits_{0}^{\min\{r,d_{\xi}\}}\sqrt{|g|}dtd\sigma
=\int\limits_{S_{x}^{n-1}}\left(\int\limits_{0}^{\min\{r,d_{\xi}\}}\det(\mathbb{A}(t,\xi))dt\right)
d\sigma,
\end{eqnarray}
where $d\sigma$ denotes the $(n-1)$-dimensional volume element on
$\mathbb{S}^{n-1}\equiv S_{x}^{n-1}\subseteq{T_{x}M}$. As in Section
1, let $r(z)=d(x,z)$ be the intrinsic distance to the point
$x\in{M}$. Since for any $\xi\in{S}_{x}^{n-1}$ and $t_{0}>0$, we
have $\nabla{r}{(\gamma_{\xi}(t_{0}))}=\gamma'_{\xi}(t_{0})$ when
the point $\gamma_{\xi}(t_{0})=\exp_{x}(t_{0}\xi)$ is away from the
cut locus of $x$ (cf. \cite{a2}), then, by the definition of a
non-zero tangent vector ``\emph{radial}" to a prescribed point on a
manifold given in the first page of \cite{KK}, we know that for
$z\in{M}\backslash(Cut(x)\cup{x})$ the unit vector field
\begin{eqnarray*}
v_{z}:=\nabla{r}{(z)}
\end{eqnarray*}
is the radial unit tangent vector at $z$. Set
\begin{eqnarray} \label{add}
l(x):=\max\limits_{z\in{M}}r(z)=\max\limits_{z\in{M}}d(x,z).
\end{eqnarray}
Then we have $l(x)=\max_{\xi}d_{\xi}$ (cf. \cite[Section 2]{fmi}).
We also need the following fact about $r(z)$ (cf. Prop. 39 on p. 266
of \cite{pp}),
\begin{eqnarray*}
\partial_{r}\Delta{r}+\frac{(\Delta{r})^2}{n-1}\leq\partial_{r}\Delta{r}+|\mathrm{Hess}r|^{2}=-\mathrm{Ric}(\partial_{r},\partial_{r}),
\qquad {\rm{with}}~~\Delta{r}=\partial_{r}\ln( \sqrt{|g|}),
\end{eqnarray*}
  with
$\partial_{r}=\nabla{r}$ as a differentiable vector (cf. Prop. 7 on
p. 47 of \cite{pp}
 for the differentiation of $\partial_{r}$), where $\Delta$ is the Laplace operator on $M$ and $\mathrm{Hess}r$ is the Hessian of $r(z)$. Then,
together with (\ref{2.2}), we have
\begin{eqnarray}
&& J''+\frac{1}{(n-1)}\mathrm{Ric}\left(\gamma'_{\xi}(t),
\gamma'_{\xi}(t)\right)J\leq 0,  \label{2.4}\\
&&J(t,\xi)=t + O(t^2), \quad  J'(t,\xi)=1+O(t). \label{2.5}
\end{eqnarray}
As shown in \cite{fmi} and also pointed out in \cite{m1}, the facts
(\ref{2.4}) and (\ref{2.5}) make a fundamental role in the
derivation of \emph{the generalized Bishop's volume comparison
theorem I} below (see Theorem \ref{theorem4} for the precise
statement). One can also find that (\ref{2.5}) is also necessary in
the proof of Theorem \ref{theorem3} in Section 3.

Denote by $inj(x)$ the injectivity radius of a point $x\in{M}$. Now,
we would like to introduce a notion of spherically symmetric
manifold which actually acts as the model space in this paper.

\begin{defn} \label{def3}
A domain $\Omega=\exp_x([0,l)\times{S}_x^{n-1}) \subset M\backslash
Cut(x)$, with $l<inj(x)$, is said to be spherically symmetric with
respect to a point $x\in \Omega$, if
 and only if
the matrix $\mathbb{A}(t,\xi)$ satisfies $\mathbb{A}(t,\xi)=h(t)I$,
for a function $h\in{C^{2}([0,l))}$, with  $h(0)=0$, $h'(0)=1$, and
$h|(0,l)>0$.
\end{defn}
Naturally, $\Omega$ in Definition \ref{def3} is a spherically
symmetric manifold and $x$ is called its \emph{base point}. Together
with (\ref{2.1}), on the set
 $\Omega$ given in Definition \ref{def3} the Riemannian metric of $M$ can be
 expressed by
 \begin{eqnarray}  \label{2.6}
 ds^{2}(\exp_{x}(t\xi))=dt^{2}+h^{2}(t)|d\xi|^{2}, \qquad
 \xi\in{S_{x}^{n-1}}, \quad 0\leq{t}<l,
 \end{eqnarray}
 with $|d\xi|^{2}$ the round metric on $\mathbb{S}^{n-1}$.
 Spherically symmetric manifolds were named as \emph{generalized space
forms} by Katz and Kondo \cite{KK}, and a standard model for such
manifolds is given by the warped product $[0,l)\times_{h}
\mathbb{S}^{n-1}$ equipped with the metric (\ref{2.6}), where $h$ is
called \emph{the warping function} and satisfies the conditions of
Definition \ref{def2}.

For a spherically symmetric manifold
$M^{\ast}:=[0,l)\times_{h}\mathbb{S}^{n-1}$ (with the base point
$p^{\ast}$) and $r<l$, by (\ref{2.3}) we have
\begin{eqnarray} \label{2.7}
\mathrm{vol}[\widetilde{B}(p^{\ast},r)]=w_{n}\int\limits_{0}^{r}h^{n-1}(t)dt,
\end{eqnarray}
and moreover, by the co-area formula (see, for instance, \cite[pp.
85-86]{ic}), we also know that the volume of the boundary
$\partial\widetilde{B}(p^{\ast},r)$ is given by
$\mathrm{vol}[\partial\widetilde{B}(p^{\ast},r)]=w_{n}h^{n-1}(r)$,
where $w_{n}$ denotes the $(n-1)$-volume of the unit sphere in
$\mathbb{R}^{n}$.

For more information about the spherically symmetric manifold
$M^{\ast}=[0,l)\times_{h}\mathbb{S}^{n-1}$ (e.g., the regularity of
the metric of $M^{\ast}$, the asymptotically spectral properties the
first Dirichlet eigenvalues of the Laplace and $p$-Laplace operators
on $M^{\ast}$, etc.), please see \cite[Section 2]{fmi} and
\cite[Section 2]{m1} in detail.

\subsection{Volume comparison theorems for manifolds with radial curvature bounded}

As before, for the given complete manifold $M$, let $d(x,\cdot)$ be
the Riemannian distance to $x$ (on $M$). In order to state volume
comparison theorems introduced below, we need the following
concepts.

\begin{defn} \label{def1}
Given a continuous function $k:[0,l)\rightarrow \mathbb{R}$, we say
that $M$ has a radial Ricci curvature lower bound $(n-1)k$ at the
point $x$ if
\begin{eqnarray*}
\mathrm{Ric}(v_z,v_z)\geq(n-1)k(d(x,z)), \quad\forall z\in
M\backslash Cut(x)\cup \{x\} ,
\end{eqnarray*}
where $\mathrm{Ric}$ is the Ricci curvature of $M$.
\end{defn}

\begin{defn} \label{def2}
Given a continuous function $k:[0,l)\rightarrow \mathbb{R}$, we say
that $M$ has a radial sectional curvature upper bound $k$ along any
unit-speed minimizing geodesic starting from a point $x\in{M}$ if
\begin{eqnarray*}
K(v_{z},V)\leq{k(d(x,z))},  ~~\forall z\in
M\backslash\left(Cut(x)\cup\{x\}\right),
\end{eqnarray*}
where $V\perp{v_{z}}$, $V\in{S^{n-1}_{z}}\subseteq{T_{z}M}$, and
$K(v_{z},V)$ is the sectional curvature of the plane spanned by
$v_z$ and $V$.
\end{defn}

\begin{remark}
\rm{As in Subsection 2.1, in Definitions \ref{def1} and \ref{def2},
$Cut(x)$ is the cut-locus of $x$ on $M$, and
$v_{z}\in{S_{p}^{n-1}}\subseteq{T_{z}M}$ is the unit tangent vector
of the minimizing geodesic $\gamma_{x,z}$ emanating from $x$ and
joining $x$ and $z$. Clearly, $v_{z}$ is in the radial direction. In
fact, the notion of having radial curvature bound has been used by
the author in \cite{fmi,m1,m2} to investigate some problems like
eigenvalue comparisons for the Laplace and $p$-Laplace operators
(between the given complete manifold and its model manifold), the
heat kernel comparison, etc. This notion can also be found in other
literatures (see, for instance, \cite{kt,st}). Let $t:=d(x,\cdot)$,
the inequality in Definition \ref{def1} (resp., Definition
\ref{def2}) becomes $\mathrm{Ric}(v_z,v_z)\geq(n-1)k(t)$ (resp.,
$K(v_{z},V)\leq{k(t)}$) for any $z\in M\backslash Cut(x)\cup \{x\}$.
We also say that the radial Ricci (resp., sectional) curvature of
$M$ is bounded from below (resp., above) by $(n-1)k(t)$ (resp.,
$k(t)$) w.r.t. $x\in{M}$ if the above inequality is satisfied. }
\end{remark}

Define a function $\widetilde{\theta}(t,\xi)$ on
$M\backslash{Cut(x)}$ as follows
\begin{eqnarray*} \label{functiontheta}
\widetilde{\theta}(t,\xi)=\left[\frac{J(t,\xi)}{h(t)}\right]^{n-1}.
\end{eqnarray*}
Then we have the following volume comparison result, which
corresponds to \cite[Theorem 3.3 and Corollary 3.4]{fmi}
(equivalently, \cite[Theorem 2.6]{m1} or \cite[Theorem 2.2.3 and
Corollary 2.2.4]{m2}).

\begin{theorem} \label{theorem4}
(A generalized Bishop's volume comparison theorem I) Given $\xi \in
S_x^{n-1}\subseteq{T_{x}M}$, and a model space
$M^{\ast}=[0,l)\times_h \mathbb{S}^{n-1}$ with the base point
$p^{\ast}$, under the curvature assumption on the radial Ricci
tensor, $\mathrm{Ric}(v_z,v_z)\geq -(n-1)h''(t)/h(t)$ on $M$, for
$z=\gamma_{\xi}(t)=\exp_{x}(t\xi)$
 with $t<\min\{d_{\xi}, l\}$,  the function
$\widetilde{\theta}$ is non-increasing in $t$. In particular, for
all $t<\min\{d_\xi, l\}$ we have $J(t,\xi)\leqslant h(t)$.
Furthermore, this inequality is strict for all $t\in (t_0,t_1]$,
with $0\leq t_0<t_1<\min \{d_{\xi},l\}$,
 if the above  curvature
assumption holds with
 a strict inequality for $t$ in the same interval. Besides, for $r_{0}<\min\{l(x),l\}$ with $l(x)$ defined by (\ref{add}), we have
$$\mathrm{vol}[B(x,r_0)]\leqslant \mathrm{vol}[\widetilde{B}(p^{\ast},r_0)],$$
with equality if and only if $B(x,r_0)$ is isometric to
$\widetilde{B}(p^{\ast},r_0)$.
\end{theorem}

Similarly, we have the following volume comparison conclusion, which
corresponds to \cite[Theorem 4.2]{fmi} (equivalently, \cite[Theorem
2.7]{m1} or \cite[Theorem 2.3.2]{m2}).

\begin{theorem}\label{theorem5}
(A generalized Bishop's volume comparison theorem II) Assume $M$ has
a radial sectional curvature upper bound $k(t)=-\frac{h''(t)}{h(t)}$
w.r.t. $x\in{M}$ for $t<\beta\leq\min\{inj_{c}(x),l\}$, where
$inj_{c}(x)=\inf_{\xi}c_{\xi}$, with  $\gamma_{\xi}(c_{\xi})$  a
first conjugate point along the geodesic
$\gamma_{\xi}(t)=\exp_{x}(t\xi)$. Then on $(0,\beta)$
\begin{eqnarray*}
\left(\frac{\sqrt{|g|}}{h^{n-1}}\right)'\geq0, \quad\quad
\sqrt{|g|}(t)\geq{h^{n-1}(t)},
 \end{eqnarray*}
and  equality occurs in the first inequality at $t_{0}\in(0,\beta)$
if and only if
 \begin{eqnarray*}
\mathcal{R}=-\frac{h''(t)}{h(t)}, \quad \mathbb{A}=h(t)I,
 \end{eqnarray*}
 on all of $[0,t_{0}]$.
\end{theorem}

\subsection{A volume comparison theorem for smooth metric measure spaces with weighted Ricci curvature bounded from below}

As mentioned at the beginning of this section, the following volume
comparison theorem proven by Wei and Wylie (cf. \cite[Theorem
1.2]{ww}) is the key point to prove Theorem \ref{theorem3}.

\begin{theorem} \cite{ww} \label{theorem6}
Let $(M,g,e^{-f}dv_{g})$ be $n$-dimensional ($n\geqslant2$) complete
smooth metric measure space with $\mathrm{Ric}_{f}\geqslant(n-1)H$.
Fix $x_{0}\in{M}$. If $\partial_{t}f\geqslant-a$ along all minimal
geodesic segments from $x_{0}$ then for $R\geqslant{r}>0$ (assume
$R\leqslant\pi/2\sqrt{H}$ if $H>0$),
\begin{eqnarray*}
\frac{\mathrm{vol}_{f}[B(x_{0},R)]}{\mathrm{vol}_{f}[B(x_{0},r)]}\leqslant{e^{aR}}\frac{\mathrm{vol}_{H}^{n}(R)}{\mathrm{vol}_{H}^{n}(r)},
\end{eqnarray*}
where $\mathrm{vol}_{H}^{n}(\cdot)$ is the volume of the geodesic
ball with the prescribed radius in the space $n$-form with constant
sectional curvature $H$, and, as before, $\mathrm{vol}_{f}(\cdot)$
denotes the weighted (or $f$-)volume of the given geodesic ball on
$M$. Moreover, equality in the above inequality holds if and only if
the radial sectional curvatures are equal to $H$ and
$\partial_{t}f\equiv-a$. In particular, if $\partial_{t}f\geqslant0$
and $\mathrm{Ric}\geqslant0$, then $M$ has $f$-volume growth of
degree at most $n$.
\end{theorem}

Therefore, given a complete and \emph{non-compact} smooth metric
measure $n$-space $(M,g,e^{-f}dv_{g})$, if $\partial_{t}f\geqslant0$
(along all minimal geodesic segments from $x_{0}$) and
$\mathrm{Ric}_{f}\geqslant0$, then by Theorem \ref{theorem6} we have
\begin{eqnarray*}
\frac{\mathrm{vol}_{f}[B(x_{0},R)]}{\mathrm{vol}_{f}[B(x_{0},r)]}\leqslant{e^{0\cdot{R}}}\cdot\frac{V_{0}(R)}{V_{0}(r)}=\frac{V_{0}(R)}{V_{0}(r)},
\end{eqnarray*}
with, as before, $V_{0}(\cdot)$ denotes the volume of the ball with
the prescribed radius in $\mathbb{R}^{n}$, which is equivalent with
\begin{eqnarray}  \label{2.9}
\frac{\mathrm{vol}_{f}[B(x_{0},R)]}{V_{0}(R)}\leqslant\frac{\mathrm{vol}_{f}[B(x_{0},r)]}{V_{0}(r)}
\end{eqnarray}
for $R\geqslant{r}>0$. Letting $r\rightarrow0$ on the right hand
side of the above inequality, and together with (\ref{2.2}),
(\ref{2.3}) and (\ref{2.5}), we can get
\begin{eqnarray*}
\frac{\mathrm{vol}_{f}[B(x_{0},R)]}{V_{0}(R)}&\leqslant&\lim\limits_{r\rightarrow0}\frac{\int\limits_{\mathbb{S}^{n-1}}\left(\int\limits_{0}^{\min\{R,d_{\xi}\}}J^{n-1}(t,\xi)\cdot{e}^{-f}dt\right)d\sigma}
{\int\limits_{\mathbb{S}^{n-1}}\int\limits_{0}^{R}t^{n-1}dtd\sigma}\nonumber\\
&=&\frac{J'(0,\xi)\cdot{e}^{-f(x_{0})}}{1}=e^{-f(x_{0})}
\end{eqnarray*}
by applying L'H\^{o}pital's rule $n$-times. Hence, if
$\partial_{t}f\geqslant0$ and $\mathrm{Ric}_{f}\geqslant0$, we have
\begin{eqnarray} \label{2.10}
\mathrm{vol}_{f}[B(x_{0},R)]\leqslant{e^{-f(x_{0})}}\cdot{V_{0}(R)}
\end{eqnarray}
for $R>0$.

\section{Proofs of main results}
\renewcommand{\thesection}{\arabic{section}}
\renewcommand{\theequation}{\thesection.\arabic{equation}}
\setcounter{equation}{0} \setcounter{maintheorem}{0}

Now, by using the facts listed in Section 2 and a similar method to
that of \cite[Theorem 1.1]{dcx}, we can prove Theorem \ref{theorem3}
as follows.

\vspace {3mm}

\textbf{\emph{Proof of Theorem \ref{theorem3}}}. Let
\begin{eqnarray} \label{3.1}
y=2a-bp+2, \qquad z=\frac{(n-2a-2)p}{2a-bp+2}=\frac{2p}{p-2},
\end{eqnarray}
where $n\geqslant3$ and $a$, $b$, $p$ are constants determined by
(\ref{1.4}). Since $t=t(\cdot):=d(x_{0},\cdot)$ is a Lipschitz
continuous function from $M$ to $\mathbb{R}$, then for any
$\lambda>0$, we can define a function $F(\lambda)$ as follows
\begin{eqnarray} \label{3.2}
F(\lambda)=\frac{p-2}{p+2}\int\limits_{M}\frac{e^{-f}dv_{g}}{t^{bp}\cdot(\lambda+t^{y})^{z-1}}.
\end{eqnarray}
by applying the Fubini theorem (cf. \cite{sy}) to (\ref{3.2}), we
have
\begin{eqnarray} \label{3.3}
F(\lambda)=\frac{p-2}{p+2}\int\limits_{0}^{+\infty}\mathrm{vol}_{f}\left[x:\frac{1}{t^{bp}\cdot(\lambda+t^{y})^{z-1}}>s\right]ds.
\end{eqnarray}
Since $\partial_{t}f\geqslant0$ (along all minimal geodesic segments
from $x_{0}$) and $\mathrm{Ric}_{f}\geqslant0$, we have (\ref{2.10})
by applying Theorem \ref{theorem6}. By making variable change
\begin{eqnarray*}
s=\frac{1}{\rho^{bp}\cdot(\lambda+\rho^{y})^{z-1}}
\end{eqnarray*}
in (\ref{3.3}) and together with (\ref{2.10}), we can obtain
\begin{eqnarray} \label{3.4}
F(\lambda)&=&\frac{p-2}{p+2}\int\limits_{0}^{+\infty}\mathrm{vol}_{f}[x:t(x)<\rho]\frac{\left[bp\lambda+(bp+(z-1)y)\rho^{y}\right]}{\rho^{bp+1}(\lambda+\rho^{y})^{z}}d\rho\nonumber\\
&=&\frac{p-2}{p+2}\int\limits_{0}^{+\infty}\mathrm{vol}_{f}[B(x_{0},\rho)]\frac{\left[bp\lambda+(bp+(z-1)y)\rho^{y}\right]}{\rho^{bp+1}(\lambda+\rho^{y})^{z}}d\rho.
\end{eqnarray}
On the other hand, by (\ref{1.4}) and (\ref{3.1}), we can get
\begin{eqnarray*}
n-bp-1>-1, \qquad n-bp-1+y(1-z)<-1.
\end{eqnarray*}
Substituting the above fact into (\ref{3.4}), it is easy to know
that $0\leqslant{F(\lambda)}<+\infty$ for any $\lambda>0$. Besides,
we also have
\begin{eqnarray*}
F'(\lambda)=-\int\limits_{M}\frac{e^{-f}dv_{g}}{t^{bp}\cdot(\lambda+t^{y})^{z}}.
\end{eqnarray*}
Therefore, from the above argument, it follows that $F$ defined by
(\ref{3.2}) is differentiable. Since for every $\lambda>0$,
$\left(\lambda+t^{y}\right)^{-\frac{z}{p}}$ is a continuous function
and tends to zero as $t\rightarrow+\infty$, which implies that there
exists at least a sequence of functions $\{g_{n}(t)\}$ in
$C^{\infty}_{0}(M)$ such that
$g_{n}(t)\rightarrow\left(\lambda+t^{y}\right)^{-\frac{z}{p}}$ as
$n\rightarrow+\infty$. By the assumption (\ref{1.5}) and an
approximation procedure for the function
$\left(\lambda+t^{y}\right)^{-\frac{z}{p}}$, we have
\begin{eqnarray*}
\left[\int\limits_{M}\frac{e^{-f}dv_{g}}{t^{bp}\cdot(\lambda+t^{y})^{z}}\right]^{\frac{2}{p}}&\leqslant&\left(\frac{yzC_{2}}{p}\right)^{2}
\int\limits_{M}\frac{e^{-f}dv_{g}}{t^{2(1+a-y)}\cdot(\lambda+t^{y})^{2+\frac{2z}{p}}}\\
&=&\left(\frac{yzC_{2}}{p}\right)^{2}\int\limits_{M}\frac{e^{-f}dv_{g}}{t^{bp-y}\cdot(\lambda+t^{y})^{z}}.
\end{eqnarray*}
Let $\ell:=\left(\frac{p}{yzC_{2}}\right)^{2}$. Then the above
inequality can be rewritten as follows
\begin{eqnarray} \label{3.5}
\ell\left[-F'(\lambda)\right]^{\frac{2}{p}}\leqslant\lambda{F}'(\lambda)+\frac{p+2}{p-2}F(\lambda).
\end{eqnarray}

Consider the function $G:(0,+\infty)\rightarrow\mathbb{R}$ defined
by
\begin{eqnarray*}
G(\lambda):=\frac{p-2}{p+2}\cdot{e}^{-f(x_{0})}\cdot\int_{\mathbb{R}^{n}}\frac{dv_{\mathbb{R}^{n}}}{|x|^{bp}\cdot(\lambda+|x|^{y})^{z-1}},
\end{eqnarray*}
where, as before, $|x|$ denotes the length of the vector
$x\in\mathbb{R}^{n}$. Since, as mentioned in Section 1, when
$C=K_{a,b}$, the extremal functions in the Caffarelli-Kohn-Nirenberg
inequality (\ref{1.2}) are of the form
$u_{\lambda}:=\left(\lambda+|x|^{y}\right)^{-\frac{z}{p}}$,
$\lambda>0$, we have
\begin{eqnarray} \label{3.6}
\left[-G'(\lambda)\right]^{\frac{2}{p}}&=&\left(e^{-f(x_{0})}\cdot\int_{\mathbb{R}^{n}}\frac{dv_{\mathbb{R}^{n}}}{|x|^{bp}\cdot(\lambda+|x|^{y})^{z}}\right)^{\frac{2}{p}}\nonumber\\
&=&\left(\frac{yzK_{a,b}}{p}\right)^{2}\left(e^{-f(x_{0})}\right)^{\frac{2}{p}}\int_{\mathbb{R}^{n}}\frac{dv_{\mathbb{R}^{n}}}{|x|^{2(1+a-y)}\cdot(\lambda+|x|^{y})^{2+\frac{2z}{p}}}\nonumber\\
&=&\left(\frac{yzK_{a,b}}{p}\right)^{2}\cdot\left[\lambda{G}'(\lambda)+\frac{p+2}{p-2}G(\lambda)\right].
\end{eqnarray}
Together with the fact $G(\lambda)=G(1)\lambda^{-\frac{2}{p-2}}$, it
follows that
\begin{eqnarray} \label{3.7}
G(1)&=&\frac{p-2}{p+2}\cdot{e}^{-f(x_{0})}\cdot\int_{\mathbb{R}^{n}}\frac{dv_{\mathbb{R}^{n}}}{|x|^{bp}\cdot(1+|x|^{y})^{z-1}}\nonumber\\
&=&2^{\frac{2}{p-2}}(p-2)\left[(n-2a-2)K_{a,b}\right]^{-\frac{2p}{p-2}}.
\end{eqnarray}
Define function $H(\lambda)$, $\lambda>0$, given by
\begin{eqnarray} \label{3.8}
H(\lambda):=A\lambda^{-\frac{2}{p-2}},
\end{eqnarray}
where $A$ satisfies
\begin{eqnarray*}
A&=&2^{\frac{2}{p-2}}(p-2)\left(\frac{\ell}{p}\right)^{\frac{p}{p-2}}\\
&=&\left(\frac{K_{a,b}}{C_{2}}\right)^{\frac{2p}{p-2}}\cdot2^{\frac{2}{p-2}}(p-2)\left[(n-2a-2)K_{a,b}\right]^{-\frac{2p}{p-2}}\\
&=&\left(\frac{K_{a,b}}{C_{2}}\right)^{\frac{n}{1+a-b}}\cdot\frac{p-2}{p+2}\cdot{e}^{-f(x_{0})}\cdot\int_{\mathbb{R}^{n}}\frac{dv_{\mathbb{R}^{n}}}{|x|^{bp}\cdot(1+|x|^{y})^{z-1}}\\
&=&\left(\frac{K_{a,b}}{C_{2}}\right)^{\frac{n}{1+a-b}}\cdot{G}(1).
\end{eqnarray*}
Clearly, by (\ref{3.7}) and (\ref{3.8}), we know that
\begin{eqnarray}  \label{3.9}
H(\lambda)=\left(\frac{K_{a,b}}{C_{2}}\right)^{\frac{n}{1+a-b}}G(\lambda).
\end{eqnarray}
Combining (\ref{3.6}) and (\ref{3.9}), one can easily check that
$H(\lambda)$ satisfies the following differential equation
\begin{eqnarray} \label{3.10}
\ell\left[-H'(\lambda)\right]^{\frac{2}{p}}=\lambda{H}'(\lambda)+\frac{p+2}{p-2}H(\lambda).
\end{eqnarray}

By L'H\^{o}pital's rule, we have
\begin{eqnarray*}
\lim\limits_{\rho\rightarrow0}\frac{\mathrm{vol}_{f}[B(x_{0},\rho)]}{V_{0}(\rho)}=\lim\limits_{\rho\rightarrow0}\frac{\int\limits_{S_{x_{0}}^{n-1}}\left(\int\limits_{0}^{\min\{\rho,d_{\xi}\}}J^{n-1}(t,\xi)\cdot{e}^{-f}dt\right)
d\sigma}{w_{n}\int\limits_{0}^{\rho}t^{n-1}dt}=e^{-f(x_{0})}.
\end{eqnarray*}
 So, for a fixed small $\epsilon>0$, there exists a number
$\eta>0$ such that
$\mathrm{vol}_{f}[B(x_{0},\rho)]\geqslant(1-\epsilon)e^{-f(x_{0})}\cdot{V}_{0}(\rho)$,
$\forall{\rho}\leqslant\eta$. Together this fact with (\ref{3.4}),
we can get
\begin{eqnarray*}
F(\lambda)&\geqslant&\frac{p-2}{p+2}(1-\epsilon)\cdot{e}^{-f(x_{0})}\cdot\int\limits_{0}^{\eta}V_{0}(\rho)\frac{\left[bp\lambda+(bp+(z-1)y)\rho^{y}\right]}{\rho^{bp+1}(\lambda+\rho^{y})^{z}}d\rho\\
&=&\frac{p-2}{p+2}(1-\epsilon)\lambda^{\frac{n+bp}{y}+1-z}\cdot{e}^{-f(x_{0})}\cdot\int\limits_{0}^{\eta\big{/}\lambda^{\frac{1}{y}}}V_{0}(s)\frac{\left[bp+(bp+(z-1)y)s^{y}\right]}{s^{bp+1}(1+s^{y})^{z}}ds\\
&=&\frac{p-2}{p+2}(1-\epsilon)\lambda^{-\frac{2}{p-2}}\cdot{e}^{-f(x_{0})}\cdot\int\limits_{0}^{\eta\big{/}\lambda^{\frac{1}{y}}}V_{0}(s)\frac{\left[bp+(bp+(z-1)y)s^{y}\right]}{s^{bp+1}(1+s^{y})^{z}}ds.
\end{eqnarray*}
On the other hand, by a direct computation, we have
\begin{eqnarray}  \label{3.11}
G(\lambda)=\frac{p-2}{p+2}\lambda^{-\frac{2}{p-2}}\cdot{e}^{-f(x_{0})}\cdot\int\limits_{0}^{+\infty}V_{0}(s)\frac{\left[bp+(bp+(z-1)y)s^{y}\right]}{s^{bp+1}(1+s^{y})^{z}}ds.
\end{eqnarray}
Therefore, it is easy to observe that
\begin{eqnarray*}
\liminf\limits_{\lambda\rightarrow0}\frac{F(\lambda)}{G(\lambda)}\geqslant1-\epsilon,
\end{eqnarray*}
and from which, one can obtain
\begin{eqnarray} \label{3.12}
\liminf\limits_{\lambda\rightarrow0}\frac{F(\lambda)}{G(\lambda)}\geqslant1
\end{eqnarray}
by letting $\epsilon\rightarrow0$.

Now, we divide into two cases to prove the assertion of Theorem
\ref{theorem3} as follows.

\emph{Case (1)}: $C_{2}>K_{a,b}$.

In this case, by (\ref{3.9}) and (\ref{3.12}), it follows that
\begin{eqnarray} \label{3.13}
\liminf\limits_{\lambda\rightarrow0}\frac{F(\lambda)}{H(\lambda)}=\left(\frac{C_{2}}{K_{a,b}}\right)^{\frac{n}{1+a-b}}\liminf\limits_{\lambda\rightarrow0}\frac{F(\lambda)}{G(\lambda)}
\geqslant\left(\frac{C_{2}}{K_{a,b}}\right)^{\frac{n}{1+a-b}}>1.
\end{eqnarray}

On the other hand, we \textbf{claim} that if there exists some
$\lambda_{0}>0$ such that $F(\lambda_{0})<H(\lambda_{0})$, then we
have $F(\lambda)<H(\lambda)$, $\forall\lambda\in(0,\lambda_{0}]$. We
will prove this by contradiction. Assume that there exists some
$\widetilde{\lambda}\in(0,\lambda_{0})$ such that
$F(\widetilde{\lambda})\geqslant{H(\widetilde{\lambda})}$. Then we
can set
\begin{eqnarray*}
\lambda_{1}:=\sup\left\{\widetilde{\lambda}<\lambda_{0}|F(\widetilde{\lambda})\geqslant{H(\widetilde{\lambda})}\right\}.
\end{eqnarray*}
So, we have $0<F(\lambda)\leqslant{H(\lambda)}$ for any
$\lambda_{1}\leqslant\lambda\leqslant\lambda_{0}$. For each
$\lambda>0$, define a function
$\phi_{\lambda}:[0,+\infty)\rightarrow\mathbb{R}$ given by
\begin{eqnarray*}
\phi_{\lambda}(m)=\ell\cdot{m}^{\frac{2}{p}}+\lambda\cdot{m}.
\end{eqnarray*}
Clearly, $\phi_{\lambda}(m)$ is increasing on $[0,+\infty)$.
Therefore, together with (\ref{3.5}) and (\ref{3.10}), it is not
difficult to get
\begin{eqnarray*}
F'(\lambda)-H'(\lambda)&\geqslant&-\phi^{-1}_{\lambda}\left(\frac{p+2}{p-2}F(\lambda)\right)+\phi^{-1}_{\lambda}\left(\frac{p+2}{p-2}H(\lambda)\right)\\
&=&\phi^{-1}_{\lambda}\left(\frac{p+2}{p-2}(H(\lambda)-F(\lambda))\right)
\geqslant\phi^{-1}_{\lambda}(0) =0
\end{eqnarray*}
for any $\lambda_{1}\leqslant\lambda\leqslant\lambda_{0}$. So,
$(F-H)'(\lambda)\leqslant0$ on $[\lambda_{1},\lambda_{0}]$.
Consequently, we can obtain
\begin{eqnarray*}
0\geqslant(F-H)(\lambda_{1})\leqslant(F-H)(\lambda_{0})<0,
\end{eqnarray*}
which is clearly a contradiction. Thus the \textbf{claim} above is
true.

By (\ref{3.13}) and the above \textbf{claim}, we have
\begin{eqnarray*}
F(\lambda)\geqslant{H(\lambda)}, \quad\forall\lambda>0.
\end{eqnarray*}
Consequently, together with (\ref{3.4}), (\ref{3.9}) and
(\ref{3.11}), we get that for any $\lambda>0$, the following
inequality
\begin{eqnarray}   \label{3.14}
\int\limits_{0}^{+\infty}\left[\mathrm{vol}_{f}[B(x_{0},\rho)]-\left(\frac{K_{a,b}}{C_{2}}\right)^{\frac{n}{1+a-b}}\cdot{e}^{-f(x_{0})}\cdot{V}_{0}(\rho)\right]\cdot\frac{\left[bp\lambda+(bp+(z-1)y)\rho^{y}\right]}{\rho^{bp+1}(\lambda+\rho^{y})^{z}}d\rho\geqslant0
\end{eqnarray}
hods. Let $b=\left(\frac{K_{a,b}}{C_{2}}\right)^{\frac{n}{1+a-b}}$.
Clearly, $0<b<1$. By Theorem \ref{theorem6}, when
$\partial_{f}\geqslant0$ (along all minimal geodesic segments from
$x_{0}$) and $\mathrm{Ric}_{f}\geqslant0$, we have (\ref{2.9}) holds
for $R\geqslant{r}>0$ and (\ref{2.10}) holds for $R>0$, which
implies that
 the volume ratio
$\mathrm{vol}_{f}[B(x_{0},\rho)]/V_{0}(\rho)$ is non-increasing for
$t\in(0,+\infty)$.
 Assume now that
$
\lim\limits_{\rho\rightarrow+\infty}\frac{\mathrm{vol}_{f}[B(x_{0},\rho)]}{e^{-f(x_{0})}\cdot{V}_{0}(\rho)}=b_{0}.
$
Clearly, $b_{0}\leqslant1$. Now, \emph{in order to get the
conclusion of Theorem \ref{theorem3} in the case $C_{2}>K_{a,b}$, it
is sufficient to show that $b_{0}\geqslant{b}$}. We will prove this
fact by contradiction. Assume that $b_{0}=b-\epsilon_{0}$ for some
$\epsilon_{0}>0$. Then there exists some $N_{0}>0$ such that
\begin{eqnarray*}
\frac{\mathrm{vol}_{f}[B(x_{0},\rho)]}{e^{-f(x_{0})}\cdot{V}_{0}(\rho)}\leqslant{b}-\frac{\epsilon_{0}}{2},
\quad \forall{t}\geqslant{N_{0}}.
\end{eqnarray*}
Substituting the above inequality into (\ref{3.14}), and together
with (\ref{2.10}), we have
\begin{eqnarray*}
&&0\leqslant\int\limits_{0}^{N_{0}}\frac{\mathrm{vol}_{f}[B(x_{0},\rho)]}{e^{-f(x_{0})}\cdot{V}_{0}(\rho)}\cdot\frac{\rho^{n}\left[bp\lambda+(bp+(z-1)y)\rho^{y}\right]}{\rho^{bp+1}(\lambda+\rho^{y})^{z}}d\rho+\\
&& \qquad
\int\limits_{N_{0}}^{+\infty}\left(b-\frac{\epsilon_{0}}{2}\right)\cdot\frac{\rho^{n}\left[bp\lambda+(bp+(z-1)y)\rho^{y}\right]}{\rho^{bp+1}(\lambda+\rho^{y})^{z}}d\rho-
b\int\limits_{0}^{+\infty}\frac{\rho^{n}\left[bp\lambda+(bp+(z-1)y)\rho^{y}\right]}{\rho^{bp+1}(\lambda+\rho^{y})^{z}}d\rho\\
&& ~
\leqslant\int\limits_{0}^{N_{0}}\frac{\rho^{n}\left[bp\lambda+(bp+(z-1)y)\rho^{y}\right]}{\rho^{bp+1}(\lambda+\rho^{y})^{z}}d\rho+
\int\limits_{N_{0}}^{+\infty}\left(b-\frac{\epsilon_{0}}{2}\right)\cdot\frac{\rho^{n}\left[bp\lambda+(bp+(z-1)y)\rho^{y}\right]}{\rho^{bp+1}(\lambda+\rho^{y})^{z}}d\rho-\\
&& \qquad \qquad
b\int\limits_{0}^{+\infty}\frac{\rho^{n}\left[bp\lambda+(bp+(z-1)y)\rho^{y}\right]}{\rho^{bp+1}(\lambda+\rho^{y})^{z}}d\rho\\
&&~=\int\limits_{0}^{N_{0}}\left(1-b+\frac{\epsilon_{0}}{2}\right)\cdot\frac{\rho^{n}\left[bp\lambda+(bp+(z-1)y)\rho^{y}\right]}{\rho^{bp+1}(\lambda+\rho^{y})^{z}}d\rho-\\
&& \qquad \qquad \frac{\epsilon_{0}}{2}\int\limits_{0}^{+\infty}\frac{\rho^{n}\left[bp\lambda+(bp+(z-1)y)\rho^{y}\right]}{\rho^{bp+1}(\lambda+\rho^{y})^{z}}d\rho\\
&&~=\int\limits_{0}^{N_{0}}\left(1-b+\frac{\epsilon_{0}}{2}\right)\cdot\frac{\rho^{n}\left[bp\lambda+(bp+(z-1)y)\rho^{y}\right]}{\rho^{bp+1}(\lambda+\rho^{y})^{z}}d\rho-\frac{n\epsilon_{0}}{2w_{n}}\cdot\frac{p+2}{p-2}\cdot{e}^{f(x_{0})}\cdot{G}(\lambda)\\
&&~\leqslant\left(1-b+\frac{\epsilon_{0}}{2}\right)\lambda^{-z}\int\limits_{0}^{N_{0}}\left[bp\lambda\rho^{n-bp-1}+(bp+(z-1)y)\rho^{n+y-bp-1}\right]d\rho-\\
&& \qquad \qquad \frac{n\epsilon_{0}}{2w_{n}}\cdot\frac{p+2}{p-2}\cdot{e}^{f(x_{0})}\cdot\lambda^{-\frac{2}{p-2}}\cdot{G}(1)\\
&&~=\left(1-b+\frac{\epsilon_{0}}{2}\right)\lambda^{-z}\left[\frac{\lambda{bp}N_{0}^{n-bp}}{n-bp}+\frac{(bp+(z-1)y)N_{0}^{n+y-bp}}{n+y-bp}\right]-\\
&& \qquad \qquad
\frac{n\epsilon_{0}(p+2)G(1)}{2w_{n}(p-2)}\cdot{e}^{f(x_{0})}\cdot\lambda^{-\frac{2}{p-2}}
\end{eqnarray*}
for every $\lambda>0$. Hence, for any $\lambda>0$, we have
\begin{eqnarray*}
0<\frac{n\epsilon_{0}(p+2)G(1)}{2w_{n}(p-2)}\cdot{e}^{f(x_{0})}\leqslant\left(1-b+\frac{\epsilon_{0}}{2}\right)\cdot\lambda^{\frac{2}{p-2}-z}\cdot\left[\frac{\lambda{bp}N_{0}^{n-bp}}{n-bp}+\frac{(bp+(z-1)y)N_{0}^{n+y-bp}}{n+y-bp}\right].
\end{eqnarray*}
Since $\frac{2}{p-2}-z+1<0$, one can get a contradiction by letting
$\lambda\rightarrow+\infty$ in the above inequality. So, this
completes the proof of the conclusion of Theorem \ref{theorem3} for
the case $C_{2}>K_{a,b}$.

\emph{Case (2)}: $C_{2}=K_{a,b}$.

In this case, by the assumption (\ref{1.5}), we have for any fixed
$\gamma>0$ that
\begin{eqnarray*}
\left(\int_{M}t^{-bp}|u|^{p}\cdot{e}^{-f}dv_{g}\right)^{\frac{1}{p}}\leqslant\left(K_{a,b}+\gamma\right)\left(\int_{M}t^{-2a}|\nabla{u}|^{2}\cdot{e}^{-f}dv_{g}\right)^{\frac{1}{2}}.
\end{eqnarray*}
Then, by the same argument to \emph{Case (1)}, we can obtain
\begin{eqnarray*}
\mathrm{vol}_{f}[B(x_{0},r)]\geqslant\left(\frac{K_{a,b}}{K_{a,b}+\gamma}\right)^{\frac{n}{1+a-b}}\cdot{e}^{-f(x_{0})}\cdot{V}_{0}(r),
\quad \forall{r}>0,
\end{eqnarray*}
which, by letting $\gamma\rightarrow0$, implies
\begin{eqnarray*}
\mathrm{vol}_{f}[B(x_{0},r)]\geqslant{e}^{-f(x_{0})}\cdot{V}_{0}(r),
\quad \forall{r}>0.
\end{eqnarray*}
This completes the proof of the conclusion of Theorem \ref{theorem3}
for the case $C_{2}=K_{a,b}$.
 \hfill $\square$

 \vspace {3mm}

Now, we give the proof of Corollary \ref{corollarys1} as follows.

\vspace {2mm}

\textbf{\emph{Proof of Corollary \ref{corollarys1}}}. By Theorem
\ref{theorem3} directly, we have
\begin{eqnarray*}
\mathrm{vol}_{f}[B(x_{0},r)]\geqslant{e}^{-f(x_{0})}\cdot{V}_{0}(r),
\quad \forall{r}>0.
\end{eqnarray*}
However, from (\ref{2.10}) which is obtained by Theorem
\ref{theorem6}, we have
\begin{eqnarray*}
\mathrm{vol}_{f}[B(x_{0},r)]\leqslant{e}^{-f(x_{0})}\cdot{V}_{0}(r),
\quad \forall{r}>0.
\end{eqnarray*}
Therefore, we have
\begin{eqnarray*}
\mathrm{vol}_{f}[B(x_{0},r)]={e}^{-f(x_{0})}\cdot{V}_{0}(r), \quad
\forall{r}>0,
\end{eqnarray*}
which, together with Theorem \ref{theorem6}, implies that the
\emph{radial} sectional curvatures are equal to $0$ and
$\partial_{t}f\equiv0$. So, we know that $f$ is a constant function
with respect to $t$, i.e., $f\equiv{f}(x_{0})$. Besides, since the
radial sectional curvatures are equal to $0$, by applying Theorems
\ref{theorem4} and \ref{theorem5} \emph{simultaneously}, we have
\begin{eqnarray*}
\mathrm{vol}[B(x_{0},r)]=V_{0}(r), \quad \forall{r}>0,
\end{eqnarray*}
and $B(x_{0},r)$ is isometric to a ball of radius $r$ in
$\mathbb{R}^{n}$ for any $r>0$, which is equivalent to say that
$(M,g)$ is isometric to
$\left(\mathbb{R}^{n},g_{\mathbb{R}^{n}}\right)$. This completes the
proof of Corollary \ref{corollarys1}.
 \hfill $\square$

\section*{Acknowledgments}
\renewcommand{\thesection}{\arabic{section}}
\renewcommand{\theequation}{\thesection.\arabic{equation}}
\setcounter{equation}{0} \setcounter{maintheorem}{0}

This work was supported by the starting-up research fund (Grant No.
HIT(WH)201320) supplied by Harbin Institute of Technology (Weihai),
the project (Grant No. HIT.NSRIF.2015101) supported by Natural
Scientific Research Innovation Foundation in Harbin Institute of
Technology, and the NSF of China (Grant No. 11401131).

 \end{document}